\documentclass[a4paper,12pt]{article}

\usepackage{hyperref}
\usepackage[pdftex]{graphicx}
\usepackage[section] {placeins}
\usepackage{enumerate}
\usepackage[normalem]{ulem}
\usepackage{amsmath}
\usepackage{amsthm}
\usepackage{latexsym}
\usepackage{amsfonts}
\usepackage{amssymb}
\usepackage{mathtools}
\usepackage{longtable}

\usepackage{tikz-cd}
\usetikzlibrary{graphs}
\usetikzlibrary{arrows.meta}
\usetikzlibrary{shapes.geometric}

\usepackage{authblk}

\usepackage{latexsym,amssymb,amsfonts, amsthm, amsmath}
\usepackage[english]{babel}
\usepackage{bm}
\usepackage{mathrsfs}
\usepackage{algpseudocode}
\usepackage{algorithm}

\algnewcommand{\IIf}[1]{\State\algorithmicif\ #1\ \algorithmicthen}
\algnewcommand{\EndIIf}{\unskip\ \algorithmicend\ \algorithmicif}

\newtheorem{thm}{Theorem}[section]
\newtheorem{lem}[thm]{Lemma}

\newtheorem{prop}[thm]{Proposition}

\newtheorem*{thm*}{Theorem \ref{th:main3}}
\newtheorem*{thm*2}{Theorem \ref{th:main1}}

\theoremstyle{definition}

\usepackage{fullpage}

\makeatletter
\newcommand{\subjclass}[2][1991]{%
\let\@oldtitle\@title%
\gdef\@title{\@oldtitle\footnotetext{#1 \emph{Mathematics subject classification.} #2}}%
}
\newcommand{\keywords}[1]{%
\let\@@oldtitle\@title%
\gdef\@title{\@@oldtitle\footnotetext{\emph{Key words and phrases.} #1.}}%
}
\makeatother

\begin{document}

\title{Weak multiset sequenceability and \\ weak BHR conjecture}
\date{}
\author[1]{Simone Costa}
\affil[1]{DICATAM, Sez.~Matematica, Universit\`a degli Studi di Brescia, Via Branze~43, I~25123 Brescia, Italy}
\affil[ ]{\texttt {simone.costa@unibs.it}}


\maketitle
\begin{abstract}
A subset $S$ of a group $(G,+)$ is $t$-{\em weakly sequenceable} if there is an ordering $(y_1, \ldots, y_k)$ of its elements such that the partial sums~$s_0, s_1, \ldots, s_k$, given by $s_0 = 0$ and $s_i = \sum_{j=1}^i y_j$ for $1 \leq i \leq k$, satisfy $s_i \neq s_j$ whenever and $1 \leq |i-j|\leq t$.
In this paper, we consider the weak sequenceability problem on multisets. In particular, we are able to prove that a multiset $M=[a_1^{\lambda_1},a_2^{\lambda_2},\dots,a_n^{\lambda_n}]$ of non-identity elements of a generic group $G$ is $t$-weakly sequenceable whenever the underlying set $\{a_1,a_2,\dots,a_n\}$ is sufficiently large (with respect to $t$) and the smallest prime divisor $p$ of $|G|$ is larger than $t$.

A related question is the one posed by the Buratti, Horak, and Rosa (briefly BHR) conjecture here considered again in the weak sense. Given a multiset $M$ and a walk $W$ in $Cay[G: \pm M]$, we say that $W$ is a realization of $M$ if $\Delta(W)=\pm M$. Here we prove that a multiset $M=[a_1^{\lambda_1},a_2^{\lambda_2},\dots,a_n^{\lambda_n}]$ of non-identity elements of $G$ admits a realization $W=(w_0,\dots,w_{\ell})$ such that $w_i\neq w_j$ whenever and $1 \leq |i-j|\leq t$ assuming that $|M|=\lambda_1+\lambda_2+\dots+\lambda_n$ is sufficiently large and the smallest prime divisor $p$ of $|G|$ is larger than $t(2t+1)$.
\end{abstract}
\section{Introduction}
A subset $S$ of a group $(G, +)$ is {\em sequenceable} if there is an ordering $(y_1, \ldots, y_k)$ of its elements such that the partial sums~$s_0, s_1, \ldots, s_k$, given by $s_0 = 0$ and $s_i = \sum_{j=1}^i y_j$ for $1 \leq i \leq k$, are distinct, with the possible exception that we may have~$s_k = s_0 = 0$.
Throughout the paper, we will use the additive notation also for non-abelian groups, i.e. $(G,+_G)$, and we avoid using the subscript in $+_G$ if it is clear from the context.

Several conjectures and questions concerning the sequenceability of subsets of groups arose in the Design Theory context: indeed this problem is related to Heffter Arrays and $G$-regular Graph Decompositions, see \cite{A15,ADMS16, CMPP18, OllisSurvey,PD}. Alspach and Liversidge combined and summarized many of them into the conjecture that if a subset of an abelian group does not contain $0$ then it is sequenceable (see \cite{AL20}).

Note that, given a sequencing of a set $S$, if $s_k\not=0$, then the partial sums define a simple path $P=(s_0,s_1,\dots,s_k)$ in the Cayley graph $Cay[G:\pm S]$ such that $\Delta(P) = \pm S$ where $\Delta(P)$ is the multiset $\pm [s_1-s_0,s_2-s_1,\dots,s_k-s_{k-1}]$. Here $Cay[G:\pm S]$ is the graph whose vertex set is $G$ and whose edges are the pairs $\{x,y\}$ such that $x-y\in \pm S$. Inspired by this interpretation, in \cite{CD} the authors proposed a weakening of this concept. In particular, they want to find an ordering of a subset $S$ of $G\setminus\{0\}$ whose partial sums $s_i$ and $s_j$ are different whenever $1 \leq |i-j|\leq t$. If such an ordering exists, we say that $S$ is $t$-{\em weakly sequenceable}. This question has been also studied in \cite{CDO} and \cite{CD1}. In particular, in \cite{CD1} the authors, using a hybrid approach that combines both Ramsey theory and the probabilistic method, proved the following asymptotic result:
\begin{thm}\label{th:main1}
Let $G$ be a generic group. Then subsets of size~$k$ of~$G\setminus\{0\}$ are $t$-weakly sequenceable whenever $k$ is large enough with respect to $t$.
\end{thm}
In this paper, we consider the weak sequenceability problem on multisets. In particular, given a multiset $M$ of non-identity elements of a group $G$, we want to find an ordering of $M$ whose partial sums $s_i$ and $s_j$ are different whenever $1 \leq |i-j|\leq t$.

Modifying the approach of \cite{CD1}, here we are able to prove that a multiset $M=[a_1^{\lambda_1},a_2^{\lambda_2},\dots,a_n^{\lambda_n}]$ (here $a_i^{\lambda_i}$ means $\lambda_i$ repetitions of $a_i$) of non-identity elements of a generic group $G$ is weakly sequenceable whenever the underlying set $\{a_1,a_2,\dots,a_n\}$ is sufficiently large and the smallest prime divisor $p$ of $|G|$ is larger than $t$.

A related question is the one first posed by Buratti for groups of prime order and then generalized by Horak and Rosa in \cite{HR09}: the so called BHR conjecture (see also \cite{OPPS}, \cite{PP1} and \cite{PP2}).
Given a multiset $M$ and a simple path $P$ in $Cay[G: \pm M]$, we say that $P$ is a realization of $M$ if $\Delta(P)=\pm M$. In this case, $M$ is said to be realizable. Roughly speaking, the difference between a realization and a sequencing of a multiset $M$ is that, in a realization, we can either add an element $a_i\in M$ or subtract it while in a sequencing we only add elements. The BHR conjecture states that any multiset $M=[a_1^{\lambda_1},a_2^{\lambda_2},\dots,a_n^{\lambda_n}]$ of $\mathbb{Z}_v$ of size $|M|=\lambda_1+\lambda_2+\dots+\lambda_n=v-1$ is realizable assuming that, for any divisor $d$ of $v$, the number of multiples of $d$ in $M$ is at most $v-d$.

In this paper we consider this conjecture, again, in the weak sense: we look for a realization of $M$ through a walk (i.e. we allow repeated vertices) and not through a simple path but we require that, locally, the vertices of the walk are different. More precisely, we prove that a multiset $M$ of non-identity elements of $G$ admits a realization $W=(w_0,\dots,w_{t})$ such that $w_i\neq w_j$ whenever $1 \leq |i-j|\leq t$ assuming that $|M|$ is sufficiently large and the smallest prime divisor $p$ of $|G|$ is larger than $t(2t+1)$.

The paper is organized as follows.
In Section $2$ we revisit the probabilistic approach of \cite{CD1}, to prove our asymptotic result on multiset weak sequenceability.

Finally, in the last section, we consider the weak variant of the BHR conjecture.
\section{Multiset weak sequenceability}
The goal of this section is to prove the asymptotic result on multiset sequenceability. For this purpose, we will recall a Ramsey-Theory kind proposition:
\begin{prop}\cite{CD1}\label{Ramsey}
Let $S$ be a subset of size $k$ of a (not necessarily abelian) group $(G, +)$ and let $t$, $\ell$ be positive integers. Then, there exists a constant $k_{t,\ell}$ such that, if $k>k_{t,\ell}$, $S$ contains a subset $T$ whose size is at least $\ell$ that does not admit zero-sum subsets (with respect to any ordering) of size $g\leq t$.
\end{prop}
With a similar proof of Proposition 2.2 of \cite{CD} we can prove the following proposition.
\begin{prop}\label{fix1}
Let $M_1=[a_1^{\alpha_1},a_2^{\alpha_2}\dots,a_n^{\alpha_n}]$ and $M_2=[a_1^{\beta_1},a_2^{\beta_2}\dots,a_n^{\beta_n}]$ be multisets of non-identity elements of a group $G$ and let us assume the underlying set of $M_2$ has size larger than $2t^2$ and that the smallest prime divisor $p$ of $|G|$ is larger than $t$. Then there is a $t$-weak sequenceable multiset $M_3$ such that
$$M_1\subseteq M_3\subseteq M_1\cup M_2$$
and $|M_3\setminus M_1|\leq t^2$. \end{prop}
\proof
First of all, we consider the sequence of $\alpha_1$ occurrences of $a_1$.
Since the order of $a_1$ is at least $p>t$, $(a_1,\dots,a_1)$ is a $t$-weak sequencing of $[a_1^{\alpha_1}]$. Now we note that the system of equations
$$\mu_1a_1+\mu_2 x\not=0:\ \mu_1\leq \min(t,\beta_1),\ \mu_2\leq \min(t,\beta_2) $$
admits a solution $x\in \{a_2,\dots,a_n\}$ provided that $n-1\geq t^2+1.$
In this case, there exists $\bar{i}$, that we may assume (without loss of generality) to be $2$, so that $(a_1,\dots,a_1, a_{2},\dots,a_{2})$ is a $t$-weak sequencing of $[a_1^{\alpha_1}, a_{2}^{\alpha_2}]$.
Reiterating, we find a $t$-weak sequencing of all the occurrences of $n-t^2$ different $a_i$s of $M_1$. Thus, we may assume (without loss of generality) that $(a_1,\dots,a_1,\dots, a_{n-t^2},\dots,a_{n-t^2})$ is a $t$-weak sequencing of $[a_1^{\alpha_1}, a_2^{\alpha_2},\dots,a_{n-t^2}^{\alpha_{n-t^2}}]$.

Now we need to place the last elements of $M_1$ that are $[a_{n-t^2+1}^{\alpha_{n-t^2+1}},\dots,a_{n}^{\alpha_n}]$. In order to do it, we denote by $m_1=|[a_1^{\alpha_1}, a_2^{\alpha_2},\dots,a_{n-t^2}^{\alpha_{n-t^2}}]|$. Then we note that the system of equations
$$s_i\not=s_{m_1}+x+\mu_{n-t^2+1}a_{n-t^2+1}:\mu_1\leq \min(t,\beta_{n-t^2+1}), |m_1+1+\mu_{n-t^2+1}-i|\leq t$$
admits a solution $x\in M_2$ provided that the underlying set of $M_2$ has size larger than $t^2$.
In this case, there exists $\bar{i}$, so that $(a_1,\dots,a_1,\dots, a_{n-t^2},\dots,a_{n-t^2},a_{\bar{i}},a_{n-t^2+1},\dots,a_{n-t^2+1})$ is a $t$-weak sequencing of $[a_1^{\alpha_1}, a_2^{\alpha_2},\dots,a_{n-t^2}^{\alpha_{n-t^2}},a_{\bar{i}},a_{n-t^2+1}^{\alpha_{n-t^2+1}}]$.
Reiterating, and assuming that the underlying set of $M_2$ has size larger than $2t^2$, we find a $t$-weak sequencing of a multiset $M_3$ such that
$$M_1\subseteq M_3\subseteq M_1\cup M_2$$
and $|M_3\setminus M_1|\leq t^2$.
\endproof

Also, now with the same proof of Proposition 2.2 of \cite{CD} we can prove the following proposition.
\begin{prop}\label{fix2}
Let $M$ and $S$ be, respectively, a multiset of size $m$ and a set of size $k$ of non-identity elements of a group $G$ and let $h$ and $t$ be positive integers such that $h\leq k-(t-1)$. Then, given a $t$-weak ordering $(y_1,\dots,y_{m})$ of $M$, there exist $h$-elements of $S$ $(y_{m+1},\dots,y_{m+h})$, such that
$$(y_1,\dots,y_{m},y_{m+1},\dots,y_{m+h})$$
is a $t$-weak ordering of
$$M\cup \{y_{m+1},\dots,y_{m+h}\}.$$
\end{prop}

Then, using a probabilistic approach, we can improve Proposition \ref{fix2} to the following one:
\begin{prop}\label{fix3}
Let $M$ and $S$ be, respectively, a multiset of size $m$ and a set of size $k$ of non-identity elements of a group $G$, and let us suppose that $k$ is large enough with respect to $t$. Then, given a $t$-weak ordering $(y_1,\dots,y_{m})$ of $M$, there exists an ordering of $S$, $(y_{m+1},\dots,y_{m+k})$, such that
$$(y_1,\dots,y_{m},y_{m+1},\dots,y_{m+k})$$
is a $t$-weak ordering of
$$M\cup S.$$
\end{prop}
\proof
Let $S$ be a subset of size $k$ of $G\setminus \{0\}$ and let $(y_1,\dots,y_{m})$ be a $t$ weak ordering of a multiset $M$. According to Proposition \ref{Ramsey}, however we choose $\ell$, if we assume that $k$ is large enough (with respect to $t$ and $\ell$), there exists a subset $T$ of $S$ whose size is $\ell$ and that does not admit zero-sum subset of size $g\leq t$. Let us set $U=S\setminus T$ and let $k'=k-\ell$ be the size of $U$. According to Proposition \ref{fix2}, given $h= k'-(t-1)=k-\ell-(t-1)$, we can order $h$ elements of $U$, namely $(y_{m+1},\dots,y_{m+h})$, in such a way that, for any ${0\leq i<j\leq \min(m+h,m+i+t)}$
$$s_i=y_1+y_2+\dots+y_i\not=y_1+y_2+\dots+y_j=s_j.$$
We denote by $U'$ the set $\{y_{m+1},\dots,y_{m+h}\}$.
Now it suffices to order the remaining $t-1$ elements of $U$ and the $\ell$ elements of $T$. Here we choose, uniformly at random an ordering $y_1,y_2,\dots,y_{m+h},z_{m+h+1},\dots, z_{m+k}$ that extend the one of $M\cup U'$. Let $X$ be the random variable that represents the number of pairs $(i, j)$ such that $s_i = s_j$ with $0\leq i<j \leq m+k$ and $j - i \leq t$. We evaluate the expected value of $X$.
Because of the linearity of the expectation, we have that:
\begin{equation}\label{E1}\mathbb{E}(X)=\sum_{\substack{0\leq i<j \leq m+h\\j - i \leq t}} \mathbb{P}(s_i=s_j)+\sum_{\substack{0\leq i\leq m+h<j\leq m+k\\j - i \leq t}} \mathbb{P}(s_i=s_j)+\sum_{\substack{m+h< i<j \leq m+k\\j - i \leq t}} \mathbb{P}(s_i=s_j).\end{equation}
Due to the choice of $U'$, when $0\leq i<j \leq m+h$, the probability $\mathbb{P}(s_i=s_j)=0$. Assuming that $0\leq i\leq m+h<j\leq m+k$,
since the equation
$$s_i=y_1+y_2+\dots+y_i =y_1+y_2+\dots+y_{m+h}+z_{m+h+1}+\dots+z_{j-1}+x=s_j$$ admits only one solution, $\mathbb{P}(s_i=s_j)$ is smaller than or equal to $\frac{1}{|S\setminus (U'\cup\{z_{h+1},\dots,z_{j-1}\})|}\leq 1/\ell.$ Here the last inequality holds because $j$ is at most $m+h+t$ and hence the cardinality of $|S\setminus (U'\cup\{z_{m+h+1},\dots,z_{m+j-1}\})|$ is at least $k-h-(t-1)=\ell$.
It follows that Equation \eqref{E1} can be rewritten as
\begin{equation}\label{E2}\mathbb{E}(X)\leq 0+\frac{t^2}{\ell}+\sum_{\substack{m+h< i<j \leq m+k\\j - i \leq t}} \mathbb{P}(s_i=s_j).\end{equation}
It is easy to see that
$$\sum_{\substack{m+h< i<j \leq m+k\\j - i \leq t}} \mathbb{P}(s_i=s_j)=\sum_{\substack{m+h< i<j \leq m+k\\j - i \leq t}} \mathbb{P}(z_{i+1}+z_{i+2}+\dots+z_j=0).$$
We divide the evaluation of this probability in two cases according to whether $\{z_{i+1},$ $z_{i+2}, \dots, z_j\}\subseteq T$ or not. In the second case, since at least one of the elements $z_{i+1},$ $z_{i+2},\dots,z_j$ is not in $T$, we can assume that $z_{i+1}\not \in T$ and since $z_{i+1}\not=0$ we can assume $j>i+1$. Reasoning as before we obtain that, under these assumptions, the probability of $z_{i+1}+z_{i+2}+\dots+z_j=0$ is at most $\frac{1}{\ell}$.
Therefore we have that:
$$\mathbb{P}(z_{i+1}+z_{i+2}+\dots+z_j=0)\leq \frac{1}{\ell}\mathbb{P}(\{z_{i+1},z_{i+2},\dots,z_j\}\not\subseteq T)+$$ $$\mathbb{P}(z_{i+1},z_{i+2},\dots,z_j\in T)\mathbb{P}(z_{i+1}+z_{i+2}+\dots+z_j=0|z_{i+1},z_{i+2},\dots,z_j\in T)$$
and, since $z_{i+1}+z_{i+2}+\dots+z_j\not=0$ whenever $z_{i+1},z_{i+2},\dots,z_j\in T$,
$$\mathbb{P}(z_{i+1}+z_{i+2}+\dots+z_j=0)\leq \frac{1}{\ell}\mathbb{P}(\{z_{i+1},z_{i+2},\dots,z_j\}\not\subseteq T).$$
Since $z_{i+1},z_{i+2},\dots,z_j$ are at most $t$ elements randomly chosen in $S\setminus U'$, the probability that all of them are contained in $T$ is at least $$\frac{{|T|\choose j-i}}{{|S\setminus U'|\choose j-i}}\geq\frac{{|T|\choose t}}{{|S\setminus U'|\choose t}} = \frac{{\ell \choose t}}{{\ell+t\choose t}}\geq \frac{(\ell-(t-1))^{t}}{(\ell+t)^{t}}.$$
It follows that
$$\mathbb{P}(z_{i+1}+z_{i+2}+\dots+z_j=0)\leq \frac{1}{\ell}\left(1-\frac{(\ell-(t-1))^{t}}{(\ell+t)^{t}}\right)$$
and hence Equation \eqref{E2} becomes
\begin{equation}\label{E3}\mathbb{E}(X)\leq \frac{t^2}{\ell}+\frac{1}{\ell}\sum_{\substack{m+h< i<j \leq m+k\\j - i \leq t}} \left(1-\frac{(\ell-(t-1))^{t}}{(\ell+t)^{t}}\right).\end{equation}
Finally, since there are less than $(\ell+t)t$ pairs $(i,j)$ such that $m+h< i<j \leq m+k$ and $j - i \leq t$, from Equation \eqref{E3} we obtain
\begin{equation}\label{E4}\mathbb{E}(X)\leq \frac{t^2}{\ell}+\frac{t(\ell+t)}{\ell}\left(1-\frac{(\ell-(t-1))^{t}}{(\ell+t)^{t}}\right).\end{equation}
Here we note that the right-hand side of the last inequality goes to zero as $\ell$ goes to infinite. Therefore for $\ell$ large enough with respect to $t$, or more precisely for $\ell\geq\bar{\ell}_t$ for a suitable $\bar{\ell}_t\in \mathbb{N}$, we have that $\mathbb{E}(X)<1$. This means that if $k$ is large enough with respect to $t$ (i.e. if $k>k_{t,\bar{\ell}_t}$) there exists an ordering on which $X=0$ that is, it is possible to extend the $t$-weak ordering of $M$ with the elements of $S$.
\endproof

We are now able to prove the main result of this section:
\begin{thm}\label{thm:mainII}
Let $(G, +)$ be a group such that the smallest prime divisor $p$ of $|G|$ is larger than $t$. Then a multiset $M=[a_1^{\lambda_1},a_2^{\lambda_2},\dots,a_n^{\lambda_n}]$ of non-identity elements of $G$ is weakly sequenceable whenever the underlying set $\{a_1,a_2,\dots,a_n\}$ is sufficiently large with respect to $t$.
\end{thm}
\proof
First of all, we note that we can split $M$ into a multisets $M_1$ and a set $S$ so that:
\begin{itemize}
\item $S=\{a_1,a_2,\dots,a_n\}$;
\item $M_1=[a_1^{\alpha_1},a_2^{\alpha_2},\dots,a_n^{\alpha_n}]$ where $\alpha_i=\lambda_i-1$.
\end{itemize}
Here, since $n$ is large enough, we can assume $n>2t^2$ and hence, due to Proposition \ref{fix1}, there is a $t$-weak sequenceable multiset $M_3$ such that
$$M_1\subseteq M_3\subseteq M_1\cup S$$
and $|M_3\setminus M_1|\leq t^2$. Denoted by $U=S\setminus M_3$ and by $k:=|U|$ its size, we have that
$$k\geq |S|-t^2=n-t^2.$$
Since $n$ is large enough, we also have that $k$ is large enough with respect to $t$. Hence, due to Proposition \ref{fix3}, it is possible to extend the $t$-weak ordering of $M_3$ with the elements of $U$. Since $M=M_3\cup U$ we obtain the thesis.
\endproof
\section{Weak BHR conjecture}
The goal of this section is to prove that any multiset $M$ of non-identity elements of $G$ admits a realization $W=(w_0,\dots,w_{t})$ such that $w_i\neq w_j$ whenever $1 \leq |i-j|\leq t$ assuming that $|M|$ is sufficiently large and the smallest prime divisor $p$ of $|G|$ is larger than $t(2t+1)$.
We begin by considering multisets $M=[a_1^{\lambda_1},a_2^{\lambda_2}]$ of two elements. We begin with a technical lemma.
\begin{lem}\label{walk1}
Let $a_1$ be an invertible element in $\mathbb{Z}_v$ and let $d$ be a divisor of $v$ such that $v\geq dt(2t+1)$.
Denoted by $W$ the walk $W=(0,a_1,2a_1,\dots,ta_1)$ we have that either
$W\cap\{1,2,\dots,dt\}=\emptyset$ or $W\cap\{-1,-2,\dots,-dt\}=\emptyset$.
\end{lem}
\proof
We consider the walk $W'=(w_0,w_1,\dots,w_{2t})$ where $w_i=ia_1$. Since $W'$ divides $\mathbb{Z}_v$ in $2t+1$ parts, there exists $\bar{i}$ such that $$W'\cap\{w_{\bar{i}}+1,w_{\bar{i}}+2,\dots,w_{\bar{i}}+\lfloor v/(2t+1)\rfloor\}=\emptyset$$ or $$W'\cap\{w_{\bar{i}}-1,w_{\bar{i}}-2,\dots,w_{\bar{i}}-\lfloor v/(2t+1)\rfloor \}=\emptyset.$$
Since $v\geq dt(2t+1)$, it follows that
$$W'\cap\{w_{\bar{i}}+1,w_{\bar{i}}+2,\dots,w_{\bar{i}}+dt\}=\emptyset$$ or $$W'\cap\{w_{\bar{i}}-1,w_{\bar{i}}-2,\dots,w_{\bar{i}}-dt \}=\emptyset.$$
Here, because of the symmetry, we may assume $\bar{i}\in \{0,1,\dots,t\}$.
It follows that the walk $W''=(w''_0,w''_1,\dots,w''_{2t-\bar{i}})$ where $w''_i=ia_1=w_{i+\bar{i}}-w_{\bar{i}}$ contains $W$ and is such that
$$W''\cap\{1,2,\dots,dt\}=\emptyset$$ or $$W''\cap\{-1,-2,\dots,-dt\}=\emptyset.$$ Hence the same property holds for the walk $W$.
\endproof
\begin{prop}\label{walk2}
Let $M=[a_1^{\lambda_1},a_2^{\lambda_2}]$ be a multiset of non-identity elements of the group $G$ and let us assume that $\lambda_1\geq 2t$ and that the smallest prime divisor $p$ of $|G|$ is larger than $t(2t+1)$.
Then, there exists a realization $W=(w_0,\dots,w_{\lambda_1+\lambda_2})$ of $M$ such that
\begin{itemize}
\item[(a)] $w_i-w_{i-1}=a_1$ for all $i\in \{1,\dots,t\}$;
\item[(b)] $w_i-w_{i-1}=a_1$ for all $i\in \{\lambda_1+\lambda_2-(t-1),\dots,\lambda_1+\lambda_2\}$;
\item[(c)] $w_i\neq w_j$ whenever $1 \leq |i-j|\leq t$.
\end{itemize}
\end{prop}
\proof
If the groups $G_1$ and $G_2$ generated respectively by $a_1$ and $a_2$ are disjoint, then they have cardinality larger than $p$ and hence than $t$. It follows that the ordering $(a_1,\dots,a_1,a_2,\dots,a_2,a_1,\dots,a_1)$ where we have $\lambda_1-t$ occurrences of $a_1$ at the beginning and $t$ occurences at the end, is a $t$-weak sequencing of $M$. Then the walk of the partial sums of this ordering provides a realization of $M$ which satisfies the conditions $(a), (b)$, and $(c)$ of this proposition.

Now, let us assume that $G_1\cap G_2=G_3\not=\{0\}$. Since both $G_1$ and $G_2$ are cyclic, we may assume $G_3$ is generated by $d\in G_1$. Also, again since $G_1$ is cyclic, we can identify $G_1$ with $\mathbb{Z}_v$ where $v=|G_1|$, $d$ with a divisor of $v$ and $a_1$ with an invertible element of $\mathbb{Z}_v$. In this case we have that $|G_3|=v/d\geq p\geq t(2t+1)$ and hence $v\geq dt(2t+1)$.

Now we consider the walk $W'=(w'_0,\dots,w'_t)=(0,a_1,2a_1,\dots,ta_1)$. Because of Lemma \ref{walk1} (applied at the end of $W'$) we have that either
$W'\cap\{ta_1+1,ta_1+2,\dots,ta_1+dt\}=\emptyset$ or $W'\cap\{ta_1-1,ta_1-2,\dots,ta_1-dt\}=\emptyset$. In the following discussion, we assume we are in the first case (the second one can be treated similarly).
Here we define the walk $W=(w_0,w_1,\dots, w_{t+\lambda_2},\dots, w_{\lambda_1+\lambda_2})$
where
$$w_i=\begin{cases} w'_i= ia_1\mbox{ if }i\leq t;\\
ta_1+(i-t)a_2 \mbox{ if } t<i\leq t+\lambda_2;\\
ta_1+\lambda_2a_2+(i-t-\lambda_2)a_1 \mbox{ if } t+\lambda_2<i\leq \lambda_1+\lambda_2.
\end{cases}
$$
Here, if $ta_1+\lambda_2a_2\not\in G_1$, $W$ is a realization of $M$ that satisfies the conditions $(a), (b)$ and $(c)$. Indeed, for construction, $W'\cap \{w_{t+1},\dots,w_{t+\min(\lambda_2,t)}\}=\emptyset$, and hence, for symmetry $$\{w_{t+\max(1,\lambda_2-t)},\dots,w_{t+\lambda_2}\}\cap\{w_{t+\lambda_2+1},w_{t+\lambda_2+2},\dots, w_{2t+\lambda_2}\}=\emptyset.$$ Moreover, since $\{w_{t+\lambda_2+1},w_{t+\lambda_2+2},\dots, w_{2t+\lambda_2}\}$ belongs to the coset of $G_1$ in $G$ represented by $w_{t+\lambda_1}=ta_1+\lambda_2a_2$, we have that $\{w_{t+\lambda_2+1},$ $w_{t+\lambda_2+2},\dots, w_{2t+\lambda_2}\}$ is also disjoint from $W'$.

Now, let us assume that $ta_1+\lambda_2a_2\in G_1$. This means that $\lambda_2a_2=\lambda_3d$ for some $\lambda_3\leq\lambda_2$. As before, we have that, for construction, \begin{equation}\label{one}
W'\cap \{w_{t+1},\dots,w_{t+\min(\lambda_2,t)}\}=\emptyset\end{equation} and, for symmetry,
$$\{w_{t+\max(1,\lambda_2-t)},\dots,w_{t+\lambda_2}\}\cap\{w_{t+\lambda_2+1},w_{t+\lambda_2+2},\dots, w_{2t+\lambda_2}\}=\emptyset.$$ It is left to prove that $w_i\not=w_j$ when $1 \leq |i-j|\leq t$, $i\leq t$ and $j\geq t+\lambda_2+1$.
Let us assume, by contradiction that there exist pairs $w_i,w_j$ with $w_i=w_j$. Let $\{i,j\}$ be chosen with minimal $j$ among such pairs. Note that, since $j-i\leq t$ it follows that $i\geq j-t\geq (t+\lambda_2+1)-t>1$. Therefore, from
$w_j=ta_1+\lambda_2a_2+(j-t-\lambda_2)a_1$ and $w_i=ia_1$, from $w_i=w_j$ it follows that also
$$w_{j-1}=ta_1+\lambda_2a_2+(j-t-\lambda_2)a_1-a_1=ia_1-a_1=w_{i-1}.$$
Due to the minimality of $j$, this means that $j-1=t+\lambda_2$. On the other hand $|(j-1)-(i-1)|=|j-i|\leq t$ and hence Equation \eqref{one} and the fact that the order of $a_2$ is at least $p>t$ impy that $w_{i-1}\not=w_{j-1}=a_1t+a_2\lambda_2$.
This is in contradiction with the assumption that $w_i=w_j$.
It follows that $W$ is a realization of $M$ that satisfies the conditions of the proposition.

The thesis follows because we can proceed analogously also when $W'\cap\{ta_1-1,ta_1-2,\dots,ta_1-dt\}=\emptyset$.
\endproof
We are now ready to consider the general case:
\begin{thm}\label{walk3}
Let $M=[a_1^{\lambda_1},a_2^{\lambda_2},\dots,a_n^{\lambda_n}]$ be a multiset of non-identity elements of the group $G$ and let us assume that $\lambda_1+\lambda_2+\dots+\lambda_n$ is large enough with respect to $t$ and that the smallest prime divisor $p$ of $|G|$ is larger than $t(2t+1)$.
Then, there exists a realization $W=(w_0,\dots,w_{\lambda_1+\lambda_2})$ of $M$ such that $w_i\neq w_j$ whenever $1 \leq |i-j|\leq t$.
\end{thm}
\proof
Since $\lambda_1+\lambda_2+\dots+\lambda_n$ is large enough, we have that either $n$ is large enough with respect to $t$ or that there exists $\bar{i}$ so that $\lambda_{\bar{i}}$ is larger than $2nt$. In the first case, the thesis follows from Theorem \ref{thm:mainII}. Therefore we can assume, without loss of generality, that $\lambda_1>2nt$.

In this case we can split $M$ into $n-1$ multisets $M_2,\dots, M_{n}$ so that
$M_i=[a_1^{\lambda_{1,i}},a_i^{\lambda_i}]$ and each $\lambda_{1,i}\geq 2t$. Because of Proposition \ref{walk2}, each $M_i$ admits a realization $W_i=(w_{0,i},\dots,w_{\lambda_{1,i}+\lambda_i,i})$ that satisfies the hypothesis. Now, in order to obtain a realization of $M$, we glue together these walks. More precisely, at first we consider $W_2=(w_{0,2},\dots,w_{\lambda_{1,2}+\lambda_2,2})$ and $W_3=(w_{0,3},\dots,w_{\lambda_{1,3}+\lambda_3,3})$ which give
$$W'=(w_{0,2},\dots,w_{\lambda_{1,2}+\lambda_2,2},w_{\lambda_{1,2}+\lambda_2,2}+w_{0,3},\dots,w_{\lambda_{1,2}+\lambda_2,2}+w_{\lambda_{1,3}+\lambda_3,3}).$$
Since $w_{i,2}-w_{i-1,2}=a_1$ whenever $i>t+\lambda_{2}$ and $w_{i,3}-w_{i-1,3}=a_1$ whenever $i\leq t$, we have that $W'$ is a realization of $M_2\cup M_3$ such that $w'_i\not=w'_j$ whenever $|i-j|\leq t$. Reiterating this procedure we obtain a realization $W=(w_0,\dots,w_{\lambda_1+\lambda_2+\dots+\lambda_n})$ of $M$ such that $w_i\neq w_j$ whenever $1 \leq |i-j|\leq t$.
\endproof
\section*{Acknowledgements} The author would like to thank Stefano Della Fiore for our useful discussions on this
topic. The author was partially supported by INdAM--GNSAGA.

\end{document}